

\baselineskip=14pt
\parskip=10pt

\font\eightrm=cmr8 
\font\eighttt=cmtt8
\magnification=\magstephalf

\def\1{{\overline{1}}}
\def\2{{\overline{2}}}
\parindent=0pt
\overfullrule=0in

\def\frac#1#2{{#1 \over #2}}
\bf
\centerline
{
Finite Analogs of Szemer\'edi's Theorem
}
\rm
\bigskip
\centerline{
{\it Paul RAFF}{$^1$} and {\it Doron ZEILBERGER}\footnote{$^1$}
{\eightrm  \raggedright
Department of Mathematics, Rutgers University (New Brunswick),
Hill Center-Busch Campus, 110 Frelinghuysen Rd., Piscataway,
NJ 08854-8019, USA.
{\eighttt [praff,zeilberg]  at math dot rutgers dot edu} ,
\hfill \break
{\eighttt http://www.math.rutgers.edu/\~{}[praff,zeilberg]} .
First written: July 13, 2009. This version: Oct. 20, 2009 (a few typos corrected, thanks
to Brian Nakamura).
Accompanied by the Maple package {\eighttt ENDRE}, 
as well as Mathematica and Java programs
downloadable from
\hfill\break
{\eighttt http://www.math.rutgers.edu/\~{}zeilberg/mamarim/mamarimhtml/szemeredi.html} .
\hfill\break
The work of DZ was supported in part by the USA National Science Foundation.
}
}

{\bf Szemer\'edi's Celebrated Theorem}

One of the crowning achievements of combinatorics is

{\bf Szemer\'edi's Theorem} ([S]): Given an integer $n \geq 1$ and an integer $k \geq 3$,
let $r_k(n)$ denote the size of any largest subset $S$ of $[n]:=\{ 1, 2, \dots, n \}$ for which
there are {\bf no} subsets of the form
$$
\{ i, i+d, i+2d, \dots, i+(k-1)d \} \quad (i \geq 1 \quad , \quad 1 \leq d < \infty ) \quad ,
$$
then $r_k(n)=o(n)$.

The depth and mainstreamness of this remarkable theorem is amply supported by the fact that at least four Fields medalists
(Klaus Roth, Jean Bourgain, Tim Gowers, and Terry Tao) and at least one Wolf prize winner (Hillel Furstenberg)
made significant contributions.

This article is yet another such contribution, and while it may not have the ``depth'' of the work of the above-mentioned
human luminaries, it does have {\it one} advantage over them. We ``cheat'' and use a computer.
It is true that, so far, we can only talk about {\it finite} analogs, but we do believe that the present approach
could be eventually extended to sharpen the current rather weak bounds.

More specifically, we prove:

{\bf Finite version of Szemer\'edi's Theorem}: Given an integer $n \geq 1$ and  integers $k \geq 3$, $D \geq 1$,
let $R_{k,D}(n)$ denote the size of any largest subset $S$ of $[n]:=\{ 1, 2, \dots, n \}$ for which
there are {\bf no} subsets of the form
$$
\{ i, i+d, i+2d, \dots, i+(k-1)d \} \quad (\,\, i \geq 1 \quad , \quad 1 \leq d \leq D \,\, ) \quad ,
$$
then there exists a rational number $\alpha_{k,D}=P_{k,D}/Q_{k,D}$ such that
$$
\lim_{n \rightarrow \infty} {{R_{k,D}(n)} \over {n}} \, = \, \alpha_{k,D} \quad .
$$

We have (rigorously!)
computed $\alpha_{k,D}$ for small $k$ and $D$ in the table below.

\vfill\eject

\halign{\bf # \hfil & # \hfil & # \hfil & # \hfil & # \hfil & # \hfil & # \hfil & # \hfil & # \hfil & # \hfil & # \hfil & # \hfil & # \hfil \cr
   & \bf 3 & \bf 4 & \bf 5 & \bf 6 & \bf 7 & \bf 8 & \bf 9 & \bf 10 & \bf 11 & \bf 12 & \bf 13 & \bf 14 \cr
1 & $2 \over 3$ & $3 \over 4$ & $4 \over 5$ & $5 \over 6$ & $6 \over 7$ & $7 \over 8$ & $8 \over 9$ & $9 \over 10$ & $10 \over 11$ & $11 \over 12$ & $12 \over 13$ & $13 \over 14$\cr
2 & $2 \over 3$ & $2 \over 3$ & $4 \over 5$ & $4 \over 5$ & $6 \over 7$ & $6 \over 7$ & $8 \over 9$ & $8 \over 9$ & $10 \over 11$ & $10 \over 11$ & $12 \over 13$ & $12 \over 13$\cr
3 & $4 \over 8$ & $8 \over 12$ & $4 \over 5$ & $4 \over 5$ & $6 \over 7$ & $6 \over 7$ & $6 \over 7$ & $20 \over 23$ & $10 \over 11$ & $10 \over 11$ & $12 \over 13$ & $12 \over 13$\cr
4 & $4 \over 9$  & $3 \over 5$ & $4 \over 5$ & $4 \over 5$ & $6 \over 7$ & $6 \over 7$ & $6 \over 7$ & $26 \over 30$ & $10 \over 11$ & $10 \over 11$ & $12 \over 13$ & $12 \over 13$\cr
5 & $4 \over 9$ & $4 \over 7$ & $16 \over 24$ & $22 \over 30$ & $6 \over 7$ &  &  &  &  &  &  & \cr
6 & $4 \over 9$ & $4 \over 7$ &  &  & &  &  &  &  &  &  & \cr
7 & $4 \over 9$ & $6 \over 11$ &  &  & &  &  &  &  &  &  & \cr
8 & $4 \over 9$ & $ 6 \over 11$&  &  & &  &  &  &  &  &  & \cr
9 & $4 \over 10$ & &  &  & &  &  &  &  &  &  & \cr
10 & $4 \over 11$ & &  &  & &  &  &  &  &  &  & \cr
11 & $8 \over 24$ & &  &  & &  &  &  &  &  &  & \cr
12 & $56 \over 177$ & &  &  & &  &  &  &  &  &  & \cr
13 & $6 \over 19$ & &  &  & &  &  &  &  &  &  & \cr
14 & $6 \over 19$ & &  &  & &  &  &  &  &  &  & \cr
15 & $6 \over 19$ & &  &  & &  &  &  &  &  &  & \cr
16  & $6 \over 19$ & &  &  & &  &  &  &  &  &  & \cr
17 & $6 \over 19$ & &  &  & &  &  &  &  &  &  & \cr
}

These numbers can get difficult to compute very quickly, but it can be seen, for example, that $\alpha_{k,1} = {k-1 \over k}$. It turns out that even more is true. $R_{k,D}(n)$ is a quasi-linear function of $n$,
and for $i=1, \dots, Q_{k,D}$ there exist integers $a_{k,D,i}$ between $0$ and $P_{k,D}-1$ such that
$$
R_{k,D} ( [Q_{k,D}] \cdot n +i)=[P_{k,D}] \cdot n + a_{k,D,i} \quad .
$$

Our proof is algorithmic, and we show how to find these explicit expressions using {\bf rigorous experimental mathematics}.

Note that $\alpha_{k,D}$ is a non-increasing sequence in $D$, and Szemer\'edi's theorem is equivalent to
the statement that
$$
\lim_{D \rightarrow \infty} \alpha_{k,D} \, = \, 0 \quad .
$$

{\bf A Wordy Formulation}

Every subset $S$ of $[1,n]=\{1,2,3, \dots, n \}$ corresponds to an $n$-letter word in the alphabet $\{0,1\}$
defined by $w[i]=1$ if and only if $i \in S$. $S$ has an arithmetical progression of size $k$ if
there is an  {\it Equidistant Letter Sequence} in the sense of the {\bf Bible Codes} of the word $1^k$
(i.e. $1$ repeated $k$ times).
Denoting by $2$ a place where the occupying letter may be either 0 or 1, we can say
that the $r_k(n)$ of Szemer\'edi's theorem defined above asks
to find the maximal number of $1$'s that an $n$-letter word in $\{0,1\}$ may have, that avoids
the {\it infinitely} many patterns
$$
(1 2^d)^{k-1} 1 \quad, \quad 0 \leq d < \infty .
$$

Analogously, the $R_{k,D}(n)$ of the finite-version 
Szemer\'edi's theorem defined above asks
to find the maximal number of $1$'s that an $n$-letter
 word in $\{0,1\}$ may have, that avoids
the {\it finitely} many patterns
$$
(1 2^d)^{k-1} 1 \quad, \quad ( \,\, 0 \leq d \leq D-1 \,\, ) .
$$

Define the {\it weight} of a word $w$ to be $t^{length}z^{\#\,\,of\,\,1s}$. Let $F_{k, D}(z,t)$ be the
weight-enumerator of all binary words avoiding the $D$ patterns
$(1 2^d)^{k-1} 1 \quad, \quad (0 \leq d \leq D-1)$. 
We will soon see that $F_{k,D}(z,t)$ is a rational function
in $(z,t)$.

Let's treat the more general case of an {\it arbitrary} set of {\it generalized patterns}.
But let's first define  {\it generalized pattern}.

{\bf Definition:} A {\it generalized pattern} is a word in the alphabet $\{0,1,2\}$, where
$2$ stands for ``space''.

Now let's say what it means to {\it contain} a pattern.

{\bf Definition:} A word $w=w_1w_2 \dots w_n$ in the alphabet $\{0,1\}$ contains the pattern $p=p_1p_2 \dots p_m$ if
there exists a position $i$ ($ 1 \leq i \leq n-m+1$) such that
$$
w_{i+j-1}=p_j \quad, \quad if \quad p_j \neq 2 \quad , \quad j=1, \dots , m \quad .
$$

For example, the word $011101101$ contains the pattern $12221$ (with $i=3$).

A word $w$ avoids a generalized pattern $p$ if it does {\it not} contain it.
A word $w$ avoids a set of generalized patterns $P$ if $w$ avoids all the members of $P$.

Analogous definitions can be made for an arbitrary finite alphabet,
where we can use SPACE ($\_$) instead of $2$. 
We will now digress to that general scenario,
and later specialize back to the binary case.

{\bf The General Problem}

Consider a finite alphabet $A$ together with a symbol
SPACE( to be denoted by $\_$)
not in $A$. We are interested in weight-enumerating
the set of words that avoid a set of patterns $P$, according to the weight
$$
weight(w_1w_2 \dots w_n)=x[w_1]x[w_2] \cdots x[w_n] \quad ,
$$
where $x[a]$ ($a \in A$) are {\it commuting indeterminates}.
For example, $weight(PAUL)=x[P]x[A]x[U]x[L]=x[A]x[L]x[P]x[U]$, $weight(DORON)=x[D]x[N]x[O]^2x[R]$.

Let $F$ be the weight-enumerator (sum of weights of its members, a formal power series in the
variables $\{x[a], a \in A\}$) of the set of such words (that avoid $P$), let's call it, for
reasons to become clear shortly, $S[P,\emptyset]$. A word belonging to it is either empty, or
else starts with one of the letters of our alphabet. If you chop that letter, what remains
is a shorter word in $S[P,\emptyset]$, but with {\it more} conditions, 
since it can not {\it start} with
a ``chopped pattern'' obtained by chopping-off the first letter for all 
those patterns of $P$
that happen to start with that letter or with $\_$ $\quad$ .

This motivates the following

{\bf Definition:}
Given a word or pattern $w=w_1 w_2 \dots w_n$, let $BEHEAD(w):=w_2 \dots w_n$.

For example,
$BEHEAD(DORON)=ORON$, $BEHEAD(PAUL)=AUL$, $BEHEAD(\_ \_ L \_ OVE)=\_ L \_ OVE$ .

Let $P$ be a set of patterns, and let $a$ be any letter of our alphabet $A$, then let
$$
P/a :=
\{ \,\,
BEHEAD(p) \, \,\,\vert \,\,\,
p \in P \quad and \quad ( \,\,
p_1=a \quad or \quad p_1=\_  \,\, ) \} \quad .
$$
For example, if the alphabet is $\{0,1\}$, and
$$
P=\{ 000, 0\_0\_0, 0\_\_0\_\_0, 111, 1\_1\_1, 1\_\_1\_\_1, \_\_101 \} \quad ,
$$
then
$$
P/0=\{ 00, \_0\_0, \_\_0\_\_0, \_101 \} \quad ,
$$
$$
P/1=\{ 11, \_1\_1, \_\_1\_\_1, \_101 \} \quad .
$$

So if $w$ belongs to our set $S[P,\emptyset]$ and it starts with the letter $a$, say, then
the chopped word obviously also avoids $P$ but in addition avoids $P/a$ at the very beginning.
This motivates us to make yet another

{\bf Definition:} Let $P$ and $P'$ be sets of patterns. The set $S[P,P']$ consists of all
words avoiding the patterns in $P$ and in addition avoiding the patterns $P'$ at the
very beginning.

Since every word in $S[P,P']$ must be either empty or else begin with one of the letters of our alphabet $A$,
we have the linear equation, for the weight-enumerators $F[P,P'](\{ x_a \})$,
$$
F[P,P']=1+ \sum_{a \in  A} x_a F[P, P/a \cup P'/a] \quad .
$$

If $P'$ contains an empty pattern, 
then of course we have the {\bf initial condition}
$F[P,P']=0$, since not even the empty word avoids the empty word as a factor.

Of course, we only care about $F[P,\emptyset]$, but in order to compute it, we need to set up
a system of linear equations featuring lots of $F[P,P']$ with 
many other (unwanted!) $P'$, but nevertheless
{\it finitely} many of them.
Since the different values of $P'$  that show up on
the right side always contain shorter patterns, and
eventually we get $P'$ that contain the empty pattern so that
we can use the initial condition, we get
{\it finitely many} (but possibly a very large number) of equations, and
{\it as many equations as unknowns}.
Also, since we know from the outset that a solution exists
(from the combinatorics), it
follows that the system of equations is non-singular, and by Cramer's rule
that we have a {\it rational function} in the variables
$$
\{ x[a] \quad \vert \quad a \in A \quad \} \quad.
$$

\vfill\eject

{\bf Specializing}

Going back to the Szemer\'edi scenario, we have a two-letter alphabet $\{0,1\}$ with weight
$x[0]=t,x[1]=zt$. For {\it any} set of forbidden patterns, in particular, those
that avoid arithmetical progression of size $k$ with spacings $\leq D$, the generating
function is of the form
$$
R(z,t)={{P(z,t)} \over {Q(z,t)}} \quad,
$$
where $t$ keeps track of the length
of words and $z$ keeps track of their number of $1$s.

Expanding $R(z,t)$ as a power-series of $t$, we get
$$
R(z,t)=\sum_{n=0}^{\infty} r_n(z )t^n \quad ,
$$
and $r_n(z)$ is a polynomial whose {\it degree} (in $z$) is the largest number $1$'s in
an $n$-letter word avoiding the set of generalized patterns.
By looking at the monomials of the denominator, $Q(z,t)$, and searching
for the monomial $z^it^j$ with {\it largest} ratio $r:=i/j$, we get that the
largest number of $1$'s in an $n$-letter word in $\{0,1\}$ is asymptotically $nr$, and more precisely,
we have the behavior described above for $R_{k,D}(n)$,
as a certain quasi-linear discrete function.

{\bf An Experimental-Yet-Rigorous Shortcut}

Solving a huge system of linear equations with {\it symbolic} coefficients is very time- and memory-
consuming. Restricting attention to the alphabet $\{0,1\}$, and letting
$f(P,P')(n)$ be the maximum number of 1's in an $n$-letter word that avoids the patterns in $P$ and
in addition, at the beginning, the patterns in $P'$, we get,
for $n>0$,
$$
f(P,P')(n)
\, = \, max( \quad
f(P, P/0 \cup P'/0) \, \, (n-1)
\quad , \quad
f(P, P/1 \cup P'/1) \,\, (n-1) \, + \, 1 \quad ) \quad .
$$

(Remember that {\it any} word in $\{0,1\}^n$, not just the one 
with the largest number of ones avoiding $P$ and $P'$, must start 
with either a 0 or a 1!). We ask the computer to {\it first} find the {\bf scheme}, in terms of a binary tree
where the left-child of $P'$ is $P/0 \cup P'/0$ and its right-child is $P/1 \cup P'/1$.
Then we ask the computer to {\it crank-out} lots of data, say, the first 
$500,000$ terms
(or whatever is needed), and then the computer automatically {\it guesses}
explicit expressions of the form
$$
R_{k,D} ([Q_{k,D}] \cdot n +i)=
[P_{k,D}] \cdot  n + a_{k,D,i} \quad , i=1 \dots Q_{k,D} \quad,
$$
for certain integers $P_{k,D}$, $Q_{k,D}$, and $a_{k,D,i}$.
Once guessed, the computer {\it automatically} gives a fully rigorous proof, {\it a posteriori},
by checking all the above equations, this time {\it symbolically}.
See the sample output of {\tt ENDRE} at the webpage of this article
for an example.

\vfill\eject

{\bf Supporting Software}

All this is implemented in the Maple package {\tt ENDRE}. 
A Mathematica program is also provided,
but only for the problems in the context of Szemer\'edi's Theorem. 
For efficiency's sake, a Java program is also available. 
See the webpage 
{\eighttt http://www.math.rutgers.edu/\~{}zeilberg/mamarim/mamarimhtml/szemeredi.html}
for these packages, as well as sample input and output.

{\bf Exact Enumeration}

From Sloane's point of view, it is interesting to crank-out as many terms as possible
of $R_{k,D}(n)$, both for their own sake, and also because they offer upper bounds
for $r_k(n)$. 
The interesting and efficient methods of the recent paper [GGK], 
that treats $r_3(n)$, may be useful to output more terms of
$R_{k,D}(n)$ for larger $D$, but of course our focus is completely different.
We do {\it symbol-crunching} rather than {\it number-crunching}.

The entries from the above table for $\alpha_{k,D}$, imply upper 
bounds for $r_4(n), r_5(n), \ldots$.

The Maple package {\tt ENDRE} also contains programs for the
{\it straight enumeration} of words of length $n$ avoiding a
set of generalized patterns, and for computing generating
functions, from which the exact asymptotics of the enumerating sequence
can be easily determined.

{\bf Finite Version of van der Waerden}

van der Waerden's theorem (for two colors) tells you that
$w_k(n)$, the number of $n$-letter words in the alphabet $\{0,1\}$,
that avoids the generalized patterns
$$
(1 2^d)^{k-1} 1 \quad, \quad
(0 2^d)^{k-1} 0 \quad, \quad
( \,\, 0 \leq d < \infty \,\, )  \quad
$$
is eventually $0$. It is still of interest to investigate
the finite version, $W_{k,D}(n)$,
the number of $n$-letter words in the alphabet $\{0,1\}$,
that avoids the generalized patterns
$$
(1 2^d)^{k-1} 1 \quad, \quad
(0 2^d)^{k-1} 0 \quad, \quad
( \,\, 0 \leq d \leq D-1 \,\, )  \quad .
$$

The Maple package {\tt ENDRE} can handle these problems as well. 

{\bf Pipe dreams}

For a fixed $k$, $\alpha_{k,D}$ gets harder and harder to compute as $D$ gets larger and larger,
{\bf but} we believe that a clever analysis of the max equations, might lead, one day,
to a {\it quantitative} understanding of how $\alpha_{k,D}$ decreases with $D$, that may (who knows?)
lead to an easier proof of Szemer\'edi's theorem, and more importantly, improved lower bounds
on $r_k(n)$.

What we are essentially doing is solving a system of recurrences of the form
$$
f_i(n)=max \, ( \, f_{a(i)}(n-1)+1 \, , \, f_{b(i)}(n-1) \, ) \quad,
$$
for $N$ sequences $\{f_i(n)\}$, $i=1 \dots N$. Here $a(i)$ $b(i)$ are some functions
from $[1,N]$ to $[1,N]$. It may be worthwhile to study such recurrences {\it for their own sake},
abstractly, and come up with a study of the asymptotic density as they depend on
$a(i),b(i)$. It is not hard to show that $f_i(n)$ can be modeled as
$$
R(Q n +i)=P n + c_i \quad,
$$
however, it is not necessarily true that $0 \leq c_i \leq P$. 
Regardless, hopefully we can get some general theorems, and since $a(i)$ and $b(i)$ are {\it arbitrary}, there is lots of elbow-room for induction.

Finally, we would
check that the {\it particular} $a(i)$, $b(i)$ that show
up satisfy some general conditions
that would enable us to get upper bounds on $\alpha_{k,D}$ as a function of $D$.

{\bf References}

[GGK] W. Gasarch, J. Glenn, C.P. Kruskal, {\it Finding large 3-free sets I: The small n case},
Journal of Computer and System Sciences {\bf 74} (2008),  628-655.
{\eighttt http://www.cs.loyola.edu/~jglenn/Papers/3apI.pdf}

[S] E. Szemer\'edi, {\it On sets of integers containing no k elements in arithmetic progression},
Acta Arith. {\bf 27}(1975), 199-245.

\end